\documentclass[conference]{IEEEtran}
\IEEEoverridecommandlockouts
\usepackage{cite}
\usepackage{amsmath,amssymb,amsfonts}
\usepackage{algorithmic}
\usepackage{multirow}
\usepackage{graphicx}
\usepackage{subcaption}

\usepackage[left=0.66in,right=0.66in,top=0.75in, bottom=1in]{geometry}

\def\BibTeX{{\rm B\kern-.05em{\sc i\kern-.025em b}\kern-.08em T\kern-.1667em\lower.7ex\hbox{E}\kern-.125emX}}

\usepackage[font=small,skip=0pt]{caption}
\usepackage[justification=centering]{caption}
\usepackage{textcomp}
\usepackage{xcolor}

\begin{document}

\title{
%Analyzing Power System Cybersecurity: 
On Graph Theory vs. Time-Domain Discrete-Event Simulation for Topology-Informed Assessment of Power Grid
%to Improve Power System 
%Cybersecurity
Cyber Risk
\thanks{K. A. Haque, L. Al Homoud, X. Zhuang, M. Elnour, A. Goulart, K. Davis, are with Texas A\&M University, College Station, TX. This work was supported by the US Department of Energy 
under award DE-CR0000018 and the National Science Foundation under grant 2220347.}}

\author{
  Khandaker Akramul Haque\IEEEauthorrefmark{1},
  Leen Al Homoud\IEEEauthorrefmark{1},
  Xin Zhuang\IEEEauthorrefmark{2},
  Mariam Elnour\IEEEauthorrefmark{1},
  %IEEE Member details
  Ana Goulart\IEEEauthorrefmark{2}, 
  Katherine Davis\IEEEauthorrefmark{1} \\
  \IEEEauthorblockA{%
    \IEEEauthorrefmark{1}Department of Electrical and Computer Engineering, Texas A\&M University, College Station, TX, USA \\
    \IEEEauthorrefmark{2}Department of Engineering Technology and Industrial Distribution, Texas A\&M University, College Station, TX, USA
  }

}

\maketitle
\thispagestyle{plain}
\pagestyle{plain}

\begin{abstract}

The shift toward more renewable energy sources and distributed generation in smart grids has underscored the significance of modeling and analyzing modern power systems as cyber-physical systems (CPS). This transformation has highlighted the importance of cyber and cyber-physical properties of modern power systems for their reliable operation. Graph theory emerges as a pivotal tool for understanding the complex interactions within these systems, providing a framework for representation and analysis. The challenge is vetting these graph theoretic 
%based
methods and other estimates of system behavior from mathematical models against reality. High-fidelity emulation and/or simulation can help %stakeholders 
answer this question, but the comparisons have been understudied. This paper employs graph-theoretic metrics to assess node risk and criticality in three distinct case studies, using 
%The studies employ SimPy, 
a Python-based discrete-event simulation called SimPy.
%is used to validate and analyze these metrics. 
%Results show that to improve network security assessment, combining graph theory and simulation is recommended to pinpoint distinct node types within the network effectively.
%\textcolor{red}{
Results for each case study show that combining graph theory and simulation
provides a topology-informed security assessment. These tools allow us to
identify critical network nodes and evaluate their performance and reliability under a cyber threat such as 
denial of service threats.%}

%By focusing on cyber and communication networks alongside physical infrastructure, this study aims to analyze the impact of vulnerability in power systems to Denial of Service (DoS) disturbances, enhancing our understanding of their cyber-physical properties and improving system resilience.

\end{abstract}

\begin{IEEEkeywords}
graph theory, betweenness centrality, cyber-physical simulation, communication networks, time-domain verification, denial-of-service
\end{IEEEkeywords}

\section{Introduction} 

Electric power systems form the backbone of the energy sector, one of the sixteen critical infrastructures defined by the U.S. Department of Homeland Security \cite{ref1}.
%, so their security is paramount.
%underscoring their importance.
%to the nation. 
%The challenge is that 
In addition, power systems are large-scale cyber-physical control systems that integrate digital communication and computation
%, and control 
technologies with physical sensing and control.  
%are also cyber-physical systems that integrate digital communication and control technologies with the physical components of a power system. 
This integration gives rise to 
%key 
functionalities
%aims 
%to enhance the 
that improve efficiency, reliability, and security of power generation, transmission, distribution, and consumption. In these 
%complex systems of systems,
%CPPS,
cyber-physical power systems,
%(CPPS), 
graph theory can offer
%provide a 
significant value to the modeling and analysis,
%of these systems, 
and its role 
%be able to
%has can 
%play a significant role in modeling and %analysis that 
can be expanded
%further utilized 
to
%help 
improve the
%holistic 
cyber-physical security and resilience of these systems more broadly. Graph theory offers a way to understand complex networks by modeling systems as graphs to assess resilience against disturbances like natural disasters and cyberattacks. This approach allows researchers to simulate scenarios and evaluate system responses to disruptions. Graph formalisms enhance experiment repeatability and comparability. Moreover, this can be adjusted to capture different system features in various models. Graph theory employs various metrics such as betweenness centrality, eccentricity centrality, eigenvector centrality, and edge betweenness centrality to analyze networks. Nodes with high betweenness and eigenvector centrality are crucial for network communication across different parts. Nodes with high edge betweenness centrality often act as bridges between network components. Conversely, nodes with low eccentricity centrality are critical within the network. These metrics collectively assess the vulnerability of infrastructure networks like cyber or power grids, offering valuable insights into their resilience and functionality. However, graph theory alone may not identify all critical network components. Therefore, we simulate real-world use cases and compare the results with graph theory metrics. These practical communication topologies can be generalized for future studies, aiding in a topology-informed security assessment of network functionality in real-world settings.

%The combined system is a combined cyber-physical power system (CPPS).  
% \textcolor{red}{However... (problem statement is missing); Thus... (the therefore is also missing).}
% \textcolor{blue}{AG: I suggest that the second paragraph should be the however... after we explain that graph theory is used, but it may not be enough to identify the critical components in the network. }

% SimPy presents an efficient and effective framework for network simulation, emphasizing the core elements while minimizing computational burden. The event-driven model adeptly manages network events, crafting a lean simulation environment tailored for Layer 3 network behaviors of the Open Systems Interconnection (OSI) model. Through asynchronous networking, SimPy captures the organic flow of events, crucial for mirroring the inherent unpredictability of real-world networks. The inclusion of packet generation and sink modeling facilitates the nuanced simulation of critical network dynamics such as traffic flow and congestion management \cite{mao2021mitigating}. Furthermore, the robust support for parallel execution and adaptable resource management simplifies the process of analyzing scalability in extensive networks. With granular simulation control, researchers and network engineers can dynamically explore protocol behaviors and network topologies, empowering them to fine-tune and optimize Layer 3 networking with precision and efficiency \cite{SimPy}. 

\begin{figure*}
    \centering
    \includegraphics[scale=0.5]{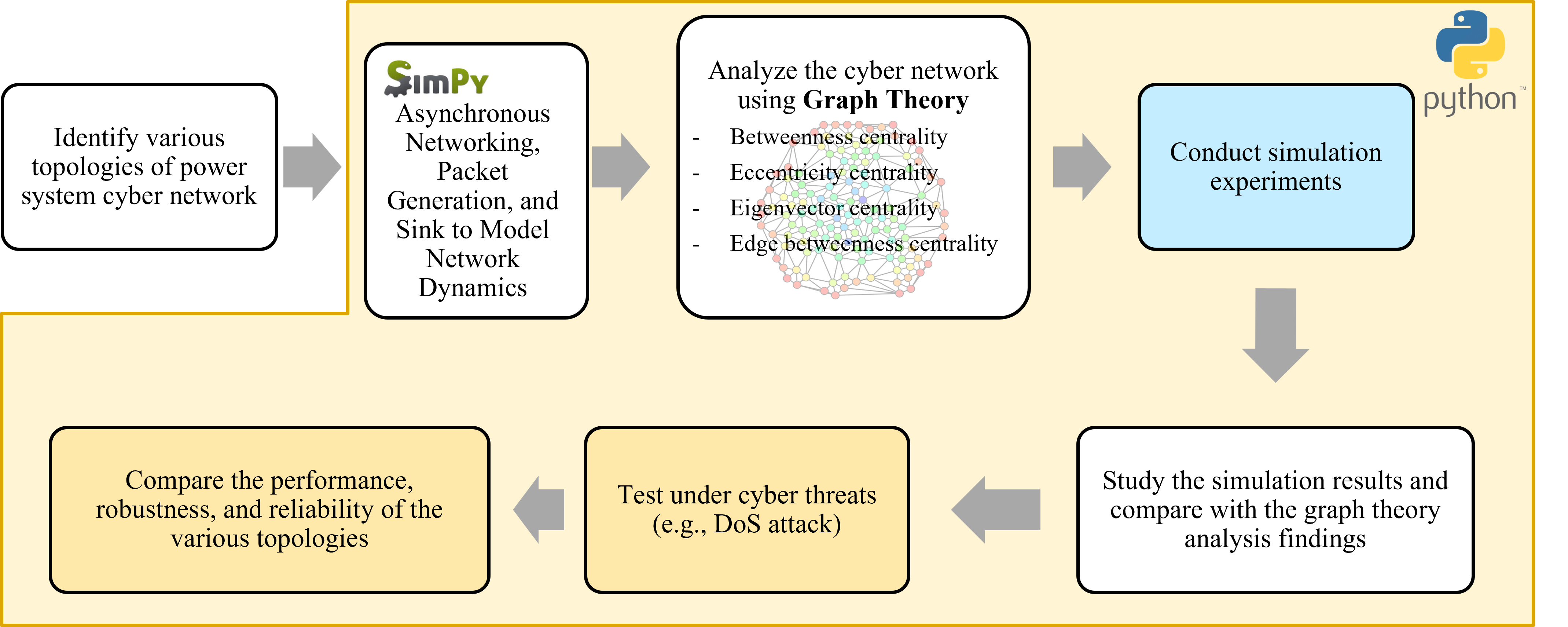}
    \caption{Methodology of power system cyber network topologies using simulation and graph theory analysis.}
    \label{fig:methodology}
\end{figure*}

This work aims to identify critical network components and ensure their functionality under anomalous conditions. We focus on evaluating the system's ability to maintain critical nodes during Denial of Service (DoS) and Distributed Denial of Service (DDoS) attacks. DoS attacks involve a single source overwhelming a target with excessive traffic, depleting resources and disrupting service. DDoS attacks, more challenging to defend against, use multiple compromised devices (botnets) to flood the target. Both types of attacks can cause service disruptions, resource exhaustion, financial losses, operational costs, reputation damage, increased vulnerability, collateral damage, and regulatory non-compliance. Using this threat model, we analyze system behavior through various case studies. Our methodology is outlined in Figure \ref{fig:methodology}.

Our simulations are based on real-world use cases, ensuring practical relevance. We use a comprehensive approach by comparing simulation results with established graph theory metrics, providing both theoretical and practical validation. The scenarios are designed for generalizability, allowing our findings to be broadly applied. We employ the reliable SimPy simulation library, ensuring methodological rigor, transparency, and reproducibility \cite{SimPy}. By integrating graph theory and simulation results, we offer a holistic and in-depth security assessment of network functionality in real-world settings. The remainder of the paper continues by discussing our motivation, contribution, and literature review in Section \ref{motivation}, methodology in Section \ref{methodology}, which details the graph metrics, the three different communication networks used for our case studies, and how the simulation is set up using SimPy in Python. The paper then explains the results and key discussions of how our findings compare across the three different case studies in Section \ref{Results}. Lastly, conclusions and future work are presented in Section \ref{Conclusion}.

\section{Motivation, Contribution, and Literature Review} \label{motivation}
This work is motivated by cyber-physical incidents targeting the Ukrainian Power Grid in 2015 and 2016 \cite{DoE}. Intruders employed spear phishing to infiltrate IT networks, subsequently manipulating control systems to induce power outages by activating breakers. This study aims to thwart such threats by leveraging cyber data to construct network graphs and simulate traffic. The goal is to pinpoint weak elements and connections within the network to bolster security and resilience. The contributions in this work are as follows:

First, we employ a process-based discrete-event simulation framework called SimPy on Python \cite{SimPy} to simulate traffic patterns and apply time-domain verification. The exponential distribution serves as the model for the behavior of individual nodes within the network in the time-domain. This distribution characterizes the inter-arrival times of packets, aligning with the fundamental property of the Poisson process. Poisson processes are widely acknowledged as effective models for aggregate traffic, particularly in scenarios involving numerous independent users exhibiting similar behavior \cite{bertsekas2021data, cleveland2000ip}. Then, we identify the vulnerable elements of the network in the event of a DoS and DDoS disturbance only in the extreme case where the adversary knows the location of the critical routers. 

Then, we compare the outcomes of the simulation experiments with those obtained using graph theory metrics, i.e.,  betweenness centrality, eccentricity centrality, eigenvector centrality and edge betweenness centrality; we evaluate the fidelity of the conducted vulnerability assessment. 

Finally, we run this analysis on three different case studies with different cyber-physical network topologies and characteristics to (i)  demonstrate  how the proposed assessment method is repeatable and generalizable and (ii) assess the cyber resilience and reliability  of the  communication network variants under this type of cyber threat. 

In today's interconnected world, securing and ensuring the reliability of power system networks is paramount. Graph theory provides valuable insights into these complex networks, aiding in calculating node vulnerability likelihood\cite{stergiopoulos2015using, umunnakwe2021cyber}. Previous research has applied graph theory for real-time vulnerability assessment in power systems and enhancing situational awareness. For instance, in \cite{biswas2020graph}, graph-based feasibility tests and network flow update schemes were employed to evaluate N-1 contingencies on meshed power networks like the IEEE 118-bus system and the Western Interconnection model. Similarly, \cite{werho2016power} used graph theory to monitor power system connectivity and identify vulnerabilities through network flow algorithms. Additionally, \cite{fan2020automated} utilized diverse graph metrics to analyze power system communication network topology, providing insights into the relationship between physical grids and communication systems. In \cite{liu2021node}, a compound centrality algorithm was introduced to assess node significance in single-layer networks, focusing on cyber-attack propagation scenarios. Moreover, \cite{adamos2023enhancing} presented a graph theory-based risk assessment strategy to enhance the resilience of cyber network architectures. Finally, \cite{de2015ranking} examined the significance of identifying central nodes in interconnected multilayer networks, emphasizing their role in processes like information propagation and congestion.

\section{Methodology}\label{methodology}

Graph theory is a valuable tool in the field of power systems for modeling, analyzing, and optimizing various aspects of electrical grids. It is used to analyze the topological structure of power systems. This includes determining the connectivity of the grid, identifying critical substations or nodes, and assessing the overall robustness of the network. According to \cite{cuzzocrea2012edge, stergiopoulos2015using}, the most used graph metrics are betweenness centrality, eccentricity centrality, eigenvector centrality and edge betweenness centrality. This section focuses on these metrics. The qualitative comparison of the three use case networks in terms of robustness and reliability is further detailed in Section~\ref{Results}. 

\subsection{Graph Metrics}
Betweenness centrality: It identifies the importance of a node in a network based on
%specifically, it %. Betweenness centrality is 
%a measure of a 
the node's centrality in %a 
the network. It quantifies how often a node acts as a bridge along the shortest path between two other nodes. Equation~\ref{eq2} shows how  to calculate betweenness centrality,
\begin{equation}
C_B(v) = \sum_{s \neq v \neq t \in V} \frac{\sigma_{st}(v)}{\sigma_{st}},
\label{eq2}
\end{equation}
 where $C_B(v)$ is the betweenness centrality of node $v$, $\sigma_{st}$ is the total number of shortest paths from node $s$ to node $t$, and $\sigma_{st}(v)$ is the number of those paths that pass through node $v$.

Eccentricity centrality: It is a metric used to assess a node's centrality within a network by measuring its distance from all other nodes. It represents the length of the longest shortest path between the node and any other node in the network, effectively indicating the node's relative position compared to others in the network structure. The eccentricity centrality, $E(v)$ of a node $v$ in network $G$ can be determined using  Equation \ref{eccentricity},
    \begin{equation}
        \label{eccentricity}
        E(v)=max_{u\in V}D(v,u)
    \end{equation}
    where $V$ is the set of all nodes in the network $G$ and $D(v,u)$ is the shortest path length between node $v$ and node $u$. The eccentricity centrality can not be calculated in a network in which every vertex is not reachable from every other vertex via directed paths. 

Eigenvector Centrality: It is a measure of centrality in a network that assigns relative scores to all nodes in the network based on the concept that connections to high-scoring nodes contribute more to the score of a node than connections to low-scoring nodes. It is a measure of influence or importance that considers both direct and indirect connections. The eigenvector centrality, $X(v)$ of a node $v$ in a network $G$ can be calculated from Equation \ref{eigenvector},
    \begin{equation}
        \label{eigenvector}
        X(v)=\frac{1}{\lambda}\sum_{u\in N(v)}X(u)
    \end{equation}
    where $N(v)$ is the set of neighbors of $v, X(u)$ is the centrality score of the neighbor node $u$ and $\lambda$ is the dominant eigenvalue of the adjacency matrix of the network.

 Edge betweenness centrality: It is used to identify the importance of edges, i.e., links that connect the nodes in a network. Identifying edges with high betweenness centrality can reveal critical components that, if disrupted, could have a significant impact. It can be calculated as in Equation~\ref{eq1},
    \begin{equation}
    \label{eq1}
    C_{B}(e) = \sum_{s \neq t \in V} \frac{\sigma_{st}(e)}{\sigma_{st}},
    \end{equation}
    where $C_{B}(e)$ is the edge betweenness centrality of edge $e$,    $\sigma_{st}$ is the total number of shortest paths from node $s$ to node $t$, $\sigma_{st}(e)$ is the number of those paths that pass through edge $e$, and $V$ is the set of nodes.

\subsection{Case Studies}
Figures \ref{fig: cyber network of IEEE 123} to \ref{fig:RingdNetwork} showcase a range of network topologies under examination. Within these networks, packet generators are depicted in red, symbolizing the origin of traffic. Routers, illustrated in blue, serve as the nodes responsible for directing traffic throughout the network. The ultimate destination of the traffic, referred to as the sink, and labeled as a balancing authority (BA), is represented in green. These network configurations were %manually extracted 
recreated from various research sources, as outlined in Table \ref{tb1}.

\begin{table}[t]
    \centering
    \scriptsize
    \caption{Summary of the case studies. 
    }
    \label{tb1}
    \begin{tabular}{cccccc} \hline
         \# & Case Study & \begin{tabular}[c]{@{}c@{}}Cyber \\Topology  \end{tabular} & Sources & Sinks & \begin{tabular}[c]{@{}c@{}}Power \\System \\Levels  \end{tabular} \\ \hline
         1 & \begin{tabular}[c]{@{}c@{}}IEEE 123 Bus Case \\ \cite{zhang2018dynamic,sahu2023reinforcement}  \end{tabular} & Mesh & Multiple & Single & \begin{tabular}[c]{@{}c@{}} Distribution \end{tabular}   \\ 
         % 2 & \begin{tabular}[c]{@{}c@{}}Wide Area Network \\in Power Systems \\ \cite{jafary2022network}  \end{tabular} & Mesh &  Single & Multiple &  \begin{tabular}[c]{@{}c@{}} Transmission \\ and \\ Distribution \end{tabular}  \\ 
         2 & \begin{tabular}[c]{@{}c@{}} Transmission and \\Distribution  Radial \\Network \cite{le2020peer} \end{tabular} & Radial &  Multiple & Single & \begin{tabular}[c]{@{}c@{}} Transmission \\ and \\ Distribution \end{tabular}   \\ 
         3 & \begin{tabular}[c]{@{}c@{}}Ring Substation \\ Topology \cite{nivethan2014modeling} \\ \end{tabular} & Ring & Multiple & Single  & \begin{tabular}[c]{@{}c@{}} Transmission \\ and \\ Distribution \end{tabular}   \\ \hline
   
    \end{tabular}
\end{table}

\subsubsection{Case 1: IEEE 123 Bus Case}

\begin{figure}[t]
\centerline{\includegraphics[scale=0.19]{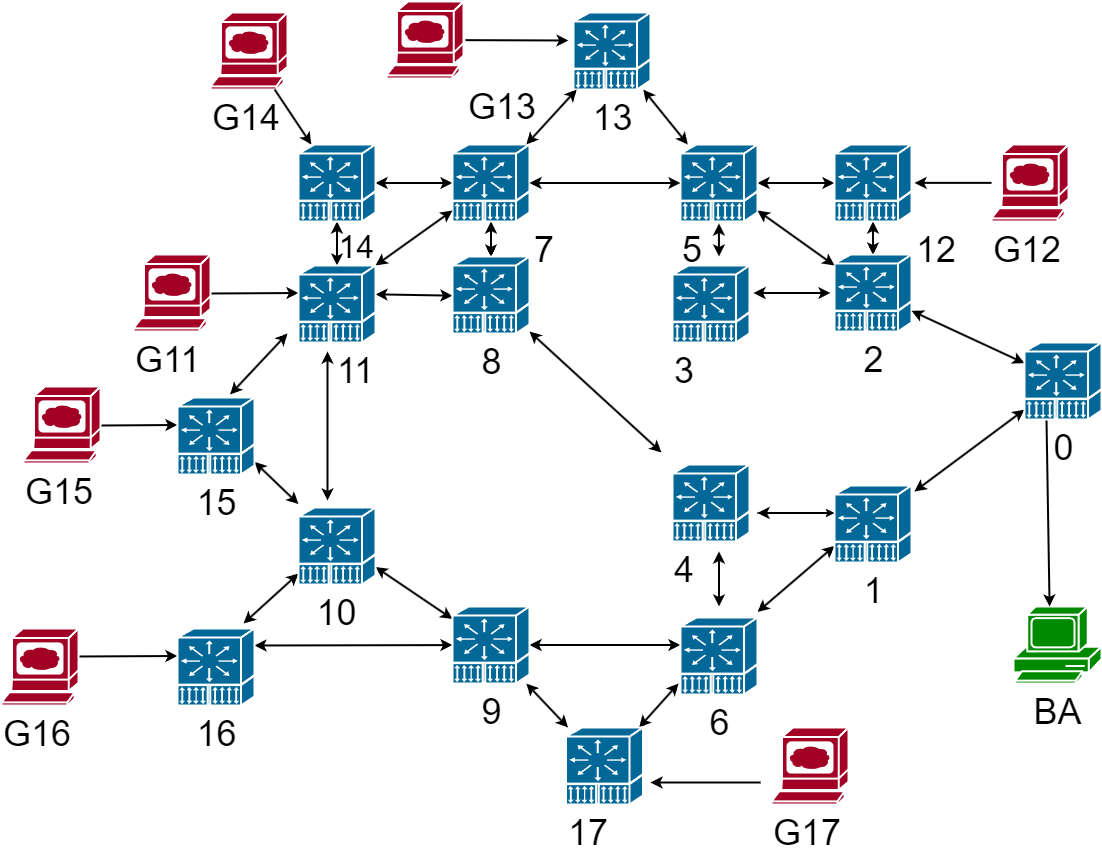}}
\caption{Case 1: The cyber network defined in~\cite{sahu2023reinforcement} based on the seven zones of the IEEE 123-bus system defined in~\cite{zhang2018dynamic}. Each source of traffic (packet generators in red) is a zone.
}
\label{fig: cyber network of IEEE 123}
\end{figure}

The IEEE 123-bus system is a widely used test system in power systems research. It represents an electrical distribution system. The communication network associated with the IEEE 123-bus system has been partitioned into seven zones for the purpose of network design and reliability analysis~\cite{zhang2018dynamic}. 
Partitioning communication networks into zones can help in managing network traffic, ensuring redundancy, and improving fault tolerance. A new network topology with these seven zones is presented in~\cite{sahu2023reinforcement}, where each zone is a source of traffic. A simplified version of the cyber network is shown in Figure \ref{fig: cyber network of IEEE 123}, where each zone consists of a packet generator (red element) that forwards the information to the centralized balancing authority acting as the sink (green element). Each packet generator is connected to a separate router (blue element), which forwards packets to the communication network. The entire communication network consists of 18 routers. 

\subsubsection{Case 2: Transmission and Distribution Radial Network}

This test case is derived from~\cite{le2020peer}, where the authors have developed a radial communication network for a combined transmission and distribution power system for market analysis studies. The root node (green node) is the interface between the transmission and
distribution network. The other nodes are in the distribution network and are ``prosumer” nodes that trade energy with their neighbors. This network topology was mapped to a power system cyber network, where the root node (green) acts as a sink, some nodes (blue) act as routers, and the rest of the edge nodes (red) act as sources that send information to the sink/BA. It is also important to note that since this case is generic, it can be overlaid and scaled for different power system models. This aids in achieving scalability for the cyber-physical analysis of large-scale power systems.

\begin{figure}[t]
    \centering    
    \includegraphics[scale = 0.20]{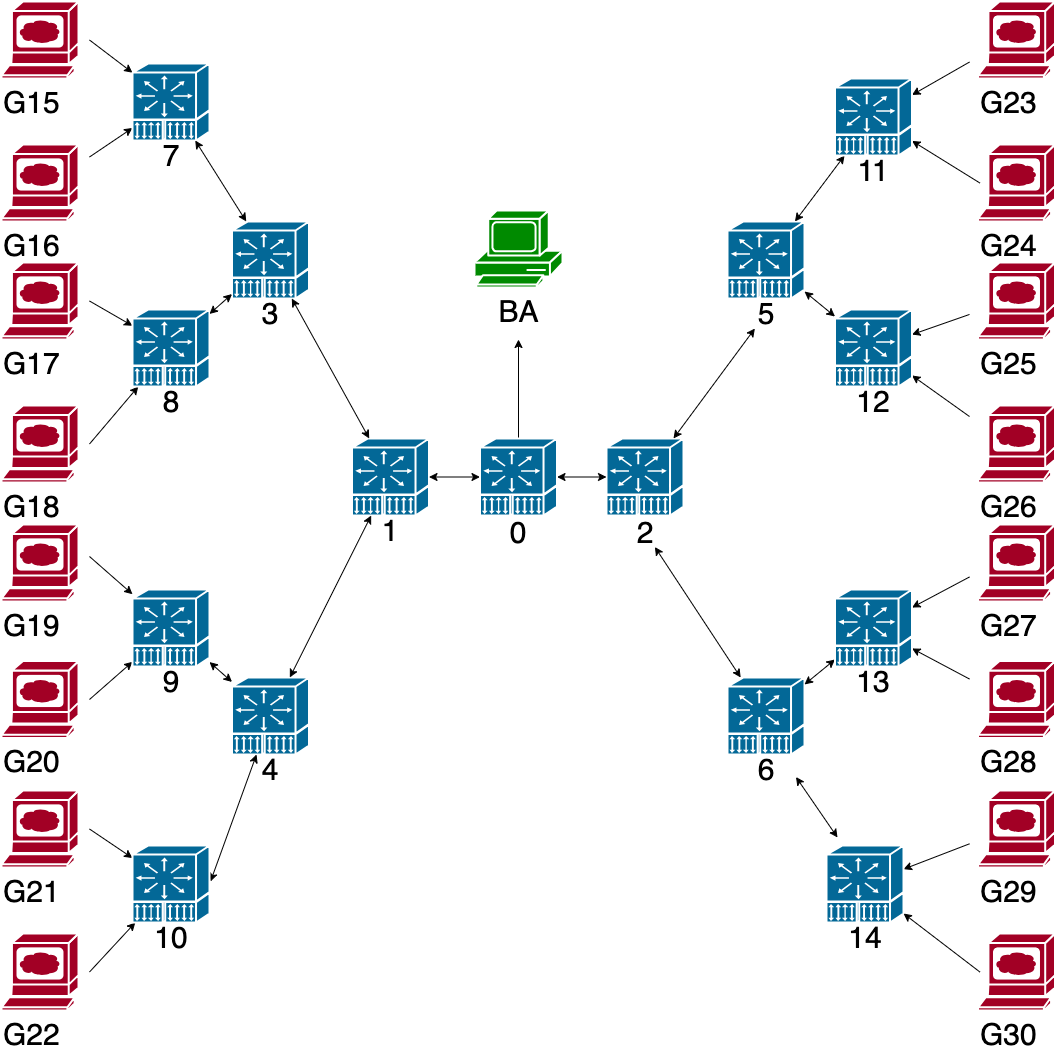}
    \caption{Case 2: A combined transmission and distribution radial communication network \cite{le2020peer}.
}
    \label{fig:radialnetwork}
\end{figure}

\subsubsection{Case 3: Ring Substation Topology} 
In a ring topology, depicted in Figure 4, devices are interconnected in a circular manner, facilitating data flow from one node to the next, typically in a single direction. This architecture is commonly utilized in wide-area networks, such as fiber optics ring topologies found in metropolitan areas. However, it is worth noting that a node losing connection can potentially disrupt the entire network.

Earlier research has indicated the utilization of ring topologies in substations adhering to the IEC 61850 standard \cite{nivethan2014modeling}. Additionally, ring configurations are applied in power distribution systems, as evidenced by \cite{ringtopologydisturbance}, where loads, wind farms, and solar farms are interconnected in a circular layout. Ring distribution systems offer heightened reliability compared to radial counterparts. Within a ring distribution setup, each load point is capable of receiving electricity from both ends of the ring, thereby enhancing system resilience.

\begin{figure}[t]
    \centering    
    \includegraphics[scale = 0.4]{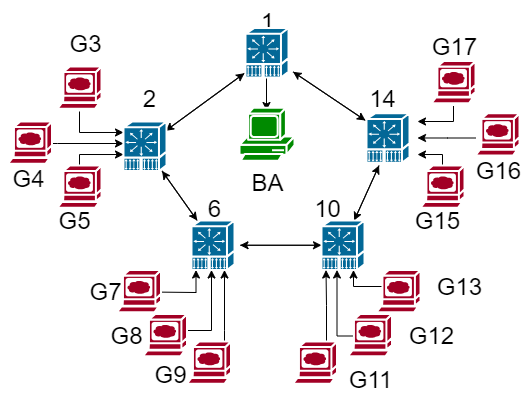}
    \caption{Case 3: A ring substation topology, similar to topology in~\cite{nivethan2014modeling}.}
    \label{fig:RingdNetwork}
\end{figure}

\subsection{Simulation Setup Overview}

To implement an event-driven simulation, a SimPy environment \cite{SimPy} is created where the generators, routers and sinks are network nodes. SimPy presents an efficient and effective framework for network simulation, emphasizing the core elements while minimizing computational burden. The event-driven model adeptly manages network events, crafting a lean simulation environment tailored for Layer 3 network behaviors of the Open Systems Interconnection (OSI) model. Through asynchronous networking, SimPy captures the organic flow of events, crucial for mirroring the inherent unpredictability of real-world networks. The inclusion of packet generation and sink modeling facilitates the nuanced simulation of critical network dynamics such as traffic flow and congestion management \cite{mao2021mitigating}. Furthermore, the robust support for parallel execution and adaptable resource management simplifies the process of analyzing scalability in extensive networks. With granular simulation control, researchers and network engineers can dynamically explore protocol behaviors and network topologies, empowering them to fine-tune and optimize Layer~3 networking with precision and efficiency. The generator in the SimPy environment generates packets based on an exponential distribution, with a packet mean size of 100 bytes. The packet inter-arrival time also follows an exponential distribution. Table~\ref{simulation_params} shows the details of our simulation setup. 

\begin{table}[h]
    \centering
    \scriptsize
    \caption{SimPy simulation parameters. 
    %\textcolor{magenta}{ Confirm if application is correct for you case, please.}
    }
    \label{simulation_params}
    \begin{tabular}{lll} \hline
         Parameter & Probability Distribution & Value  \\ \hline
         Packet Size & Exponential & Mean = 100 Bytes \\ 
         Packet inter-arrival time & Exponential & Max = 2 sec  \\ 
         Router's port rate  & Exponential & 2.2 packets/sec \\ 
         Sampling rate & Exponential & Max = 0.5 samples/sec \\ \hline
    \end{tabular}
\end{table}

\begin{table*}[t]
\centering
\caption{
Rank of critical routers and edges between routers based on $C_B(v)$, $E(v)$, $X(v)$ and $C_B(e)$.
%
%In the results, 
%($/$) denotes~\textcolor{red}{...[please help clarify…?} Also what are they being sorted on? What is "Rank" go through each column-- and clarifying the "Router" "CB/E/X" and "Edge" and "CB(e)" columns each, a little more concretely would help..] between routers represents the rank based on centrality.
%
%The 
The `Router' column gives
%indicates 
the router ranked at position `Rank,' based on the metric `CB(v)/E(v)/X(v),' where the position of the ($/$) in the router column denotes which of the
%these 
three metrics its rank is based on.
%its rank is based on.
%
For example, for Case 1, in row 1, $7/X/7$ and $0.28/X/.044$ indicate that Router 7 is ranked first based on a $C_B(v)$ metric of 0.28 and 
%ranked
first based on $X(v)$ metric of 0.44, and `X' denotes that $E(v)$ cannot be calculated since the graph is not strongly connected. `Edge' and `CB(e)' indicate the edge rankings and values, and
%the
%, using
`()'
%to 
denotes the clusters with elements of the same ranking.
%that edge 5-7 is ranked first based on CB(e) of 0.1.
%
%In Case 1, 
%
%Router 7 is first based on 0.28
%
%
%NA in case 2 - all routers have been classified
%
%Case 3:
%Case 2 and 3 are symmetrical
%
%only has 5 routers hence 
The NAs at the end
%} 
%
%In Case 3, 
%NA indicates 
denote not applicable, where all elements have already been classified.
%\textcolor{red}{why is it not applicable?}. 
%For example, 
%\textcolor{red}{[please help clarify: ]  in case 2 router 1, 2 is rank 1 for $C_B(v)$, router 0 is rank 1 for $E(v)$ and router 1, 2 is rank 1 for $X(v)$.
%}
%Rank of critical routers and edge between routers based on betweenness, $C_B(v)$/eccentricity, $E(v)$/eigenvector, $X(v)$ centrality and edge betweenness centrality, $C_B(e)$ metric where X represents the value of centrality which cannot be calculated since the graph is not strongly connected (Case 1) and NA represents not applicable; ($/$) between routers represents the rank based on centrality. For example, in case 2 router 1, 2 is rank 1 for $C_B(v)$, router 0 is rank 1 for $E(v)$ and router 1, 2 is rank 1 for $X(v)$.
}
\label{tb2}
\setlength{\arrayrulewidth}{0.6mm}
\scriptsize
\begin{tabular}{|c||cccc||cccc||cccc|}
\hline
\multirow{2}{*}{\textbf{Rank}} & \multicolumn{4}{c||}{\textbf{Case 1}}                                                                                                                                                  & \multicolumn{4}{c||}{\textbf{Case 2}}                                                                                                                                                                                                                                                                                                                                          & \multicolumn{4}{c|}{\textbf{Case 3}}                                                                                                                                                                                                               \\ \cline{2-13} 
                               & \multicolumn{1}{c|}{\textbf{Router}} & \multicolumn{1}{c|}{\textbf{\begin{tabular}[c]{@{}c@{}}CB(v)/E(v)/\\ X(v)\end{tabular}}} & \multicolumn{1}{c|}{\textbf{Edge}} & \textbf{CB(e)} & \multicolumn{1}{c|}{\textbf{Router}}                                                                                                       & \multicolumn{1}{c|}{\textbf{\begin{tabular}[c]{@{}c@{}}CB(v)/E(v)/\\ X(v)\end{tabular}}} & \multicolumn{1}{c|}{\textbf{Edge}}                                                                                   & \textbf{CB(e)} & \multicolumn{1}{c|}{\textbf{Router}}                                                              & \multicolumn{1}{c|}{\textbf{\begin{tabular}[c]{@{}c@{}}CB(v)/E(v)/\\ X(v)\end{tabular}}} & \multicolumn{1}{c|}{\textbf{Edge}} & \textbf{CB(e)} \\ \hline
1                              & \multicolumn{1}{c|}{7/X/7}           & \multicolumn{1}{c|}{\begin{tabular}[c]{@{}c@{}}0.28/X/\\ 0.44\end{tabular}}              & \multicolumn{1}{c|}{5-7}           & 0.10           & \multicolumn{1}{c|}{\begin{tabular}[c]{@{}c@{}}(1, 2)/0/\\    \\ (1, 2)\end{tabular}}                                                          & \multicolumn{1}{c|}{\begin{tabular}[c]{@{}c@{}}0.63/4/\\ 0.35\end{tabular}}              & \multicolumn{1}{c|}{(0-2, 0-1)}                                                                                        & 0.52           & \multicolumn{1}{c|}{\begin{tabular}[c]{@{}c@{}}(6, 10)/\\ (1, 2, 6, \\ 10, 14)/\\ (6, 10)\end{tabular}} & \multicolumn{1}{c|}{\begin{tabular}[c]{@{}c@{}}0.45/3/\\ 0.43\end{tabular}}              & \multicolumn{1}{c|}{(6-10)}          & 0.31           \\ \hline
2                              & \multicolumn{1}{c|}{11/X/11}         & \multicolumn{1}{c|}{\begin{tabular}[c]{@{}c@{}}0.22/X/\\ 0.42\end{tabular}}              & \multicolumn{1}{c|}{7-11}          & 0.08           & \multicolumn{1}{c|}{\begin{tabular}[c]{@{}c@{}}0/(1, 2)/ \\ (3, 4, 5, \\ 6)\end{tabular}}                                                   & \multicolumn{1}{c|}{\begin{tabular}[c]{@{}c@{}}0.52/5/\\ 0.29\end{tabular}}              & \multicolumn{1}{c|}{\begin{tabular}[c]{@{}c@{}}(2-5, 2-6,\\    \\ 1-3, 1-4)\end{tabular}}                              & 0.36           & \multicolumn{1}{c|}{\begin{tabular}[c]{@{}c@{}}(14, 2)/\\ NA/\\ (14, 2)\end{tabular}}                 & \multicolumn{1}{c|}{\begin{tabular}[c]{@{}c@{}}0.39/NA/\\ 0.38\end{tabular}}             & \multicolumn{1}{c|}{(2-6,   10-14)}  & 0.26           \\ \hline
3                              & \multicolumn{1}{c|}{5/X/5}           & \multicolumn{1}{c|}{\begin{tabular}[c]{@{}c@{}}0.20/X/\\ 0.33\end{tabular}}              & \multicolumn{1}{c|}{10-11}         & 0.07           & \multicolumn{1}{c|}{\begin{tabular}[c]{@{}c@{}}(3, 4, 5, \\ 6)/\\ (3, 4, \\ 5, 6)/\\ (7, 8, 9, \\ 10, 11, 12, \\ 13, 14)\end{tabular}}           & \multicolumn{1}{c|}{\begin{tabular}[c]{@{}c@{}}0.35/6/\\ 0.18\end{tabular}}              & \multicolumn{1}{c|}{\begin{tabular}[c]{@{}c@{}}(3-7, 3-8, \\ 4-9, 4-10, \\ 5-11, \\ 5-12, \\ 6-13, 6-14)\end{tabular}} & 0.18           & \multicolumn{1}{c|}{1/NA/1}                                                                       & \multicolumn{1}{c|}{\begin{tabular}[c]{@{}c@{}}0.23/NA/\\ 0.30\end{tabular}}             & \multicolumn{1}{c|}{(1-2,   1-14)}   & 0.20           \\ \hline
4                              & \multicolumn{1}{c|}{10/X/8}          & \multicolumn{1}{c|}{\begin{tabular}[c]{@{}c@{}}0.15/X/\\ 0.29\end{tabular}}              & \multicolumn{1}{c|}{8-4}           & 0.06           & \multicolumn{1}{c|}{\begin{tabular}[c]{@{}c@{}}(7, 8, 9, \\ 10, 11, 12, \\ 13, 14)/\\ (7, 8, 9, \\ 10, 11, 12, \\ 13, 14)/ \\ NA\end{tabular}} & \multicolumn{1}{c|}{\begin{tabular}[c]{@{}c@{}}0.13/7/\\ NA\end{tabular}}                & \multicolumn{1}{c|}{NA}                                                                                              & NA             & \multicolumn{1}{c|}{NA}                                                                           & \multicolumn{1}{c|}{NA}                                                                  & \multicolumn{1}{c|}{NA}            & NA             \\ \hline
\end{tabular}
\end{table*}

\begin{figure*}[t]
\scriptsize
\centering
\begin{subfigure}{0.31\textwidth}
    \includegraphics[width=\textwidth]{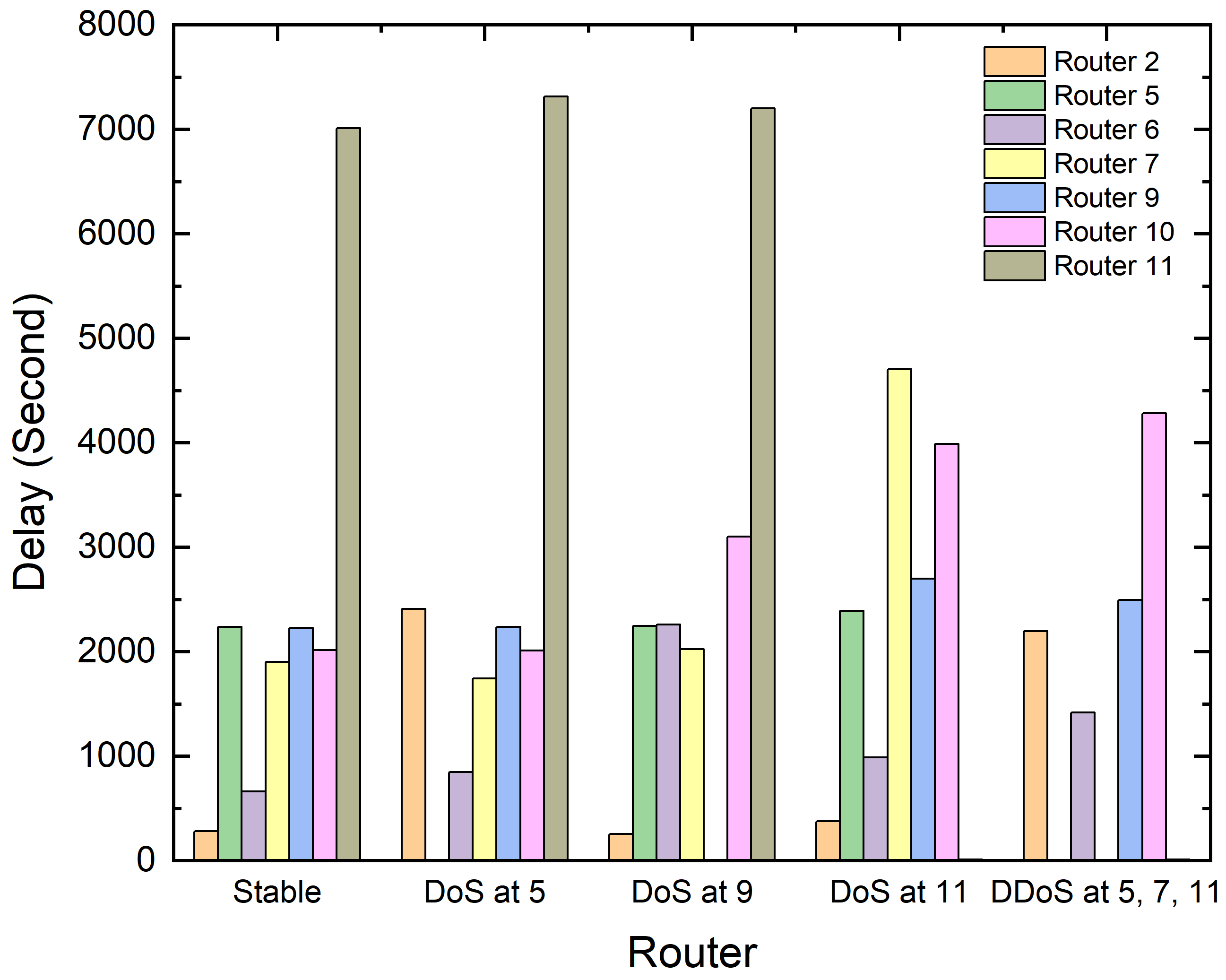}
    \caption{Case 1}
    \label{fig:first}
\end{subfigure}
\hfill
% \begin{subfigure}{0.42\textwidth}
%     \includegraphics[width=\textwidth]{figs/Case2_delay.png}
%     \caption{Case 2}
%     \label{fig:second}
% \end{subfigure}
% \hfill
\begin{subfigure}{0.31\textwidth}
    \includegraphics[width=\textwidth]{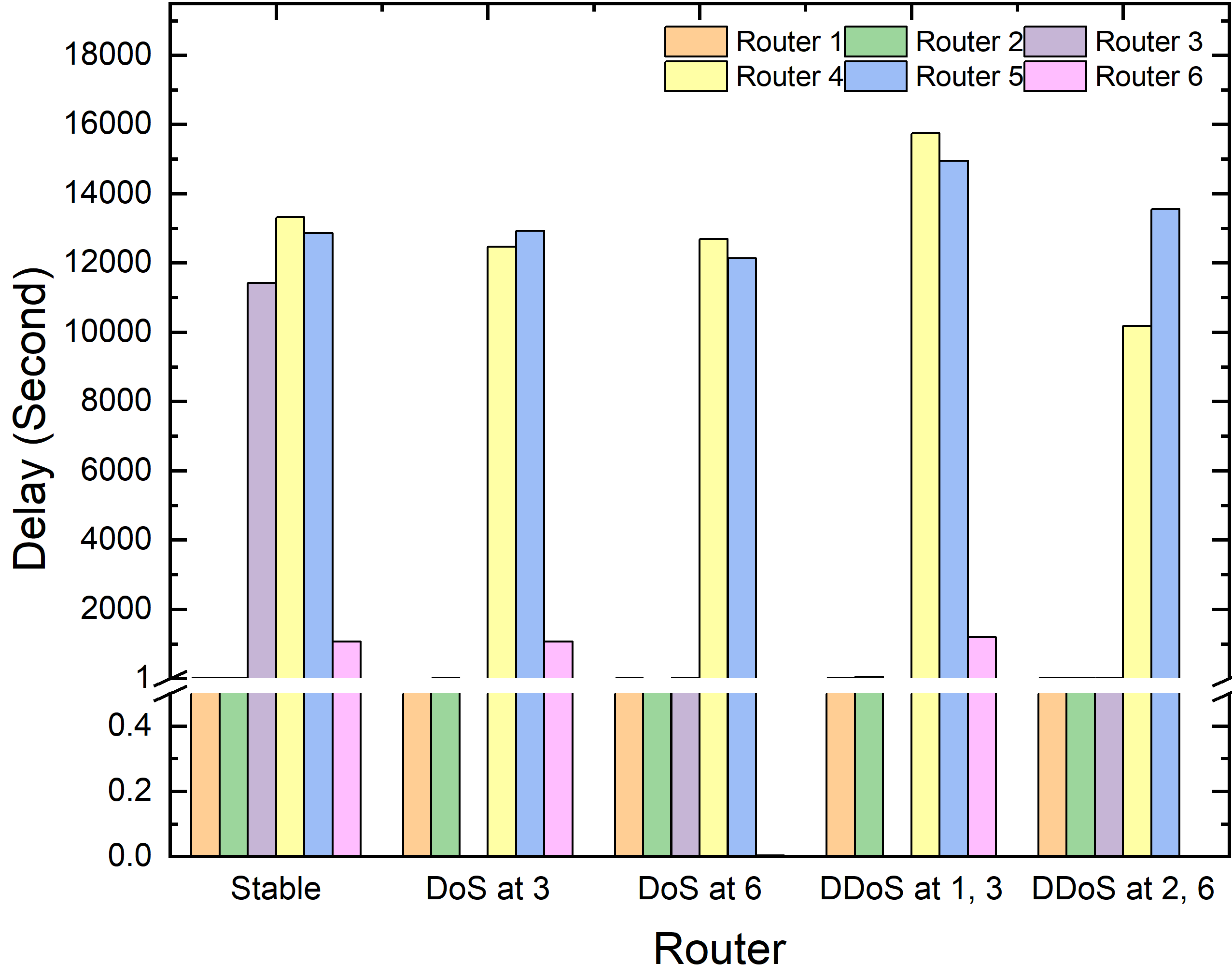}
    \caption{Case 2}
    \label{fig:third}
\end{subfigure}
\hfill
\begin{subfigure}{0.31\textwidth}
    \includegraphics[width=\textwidth]{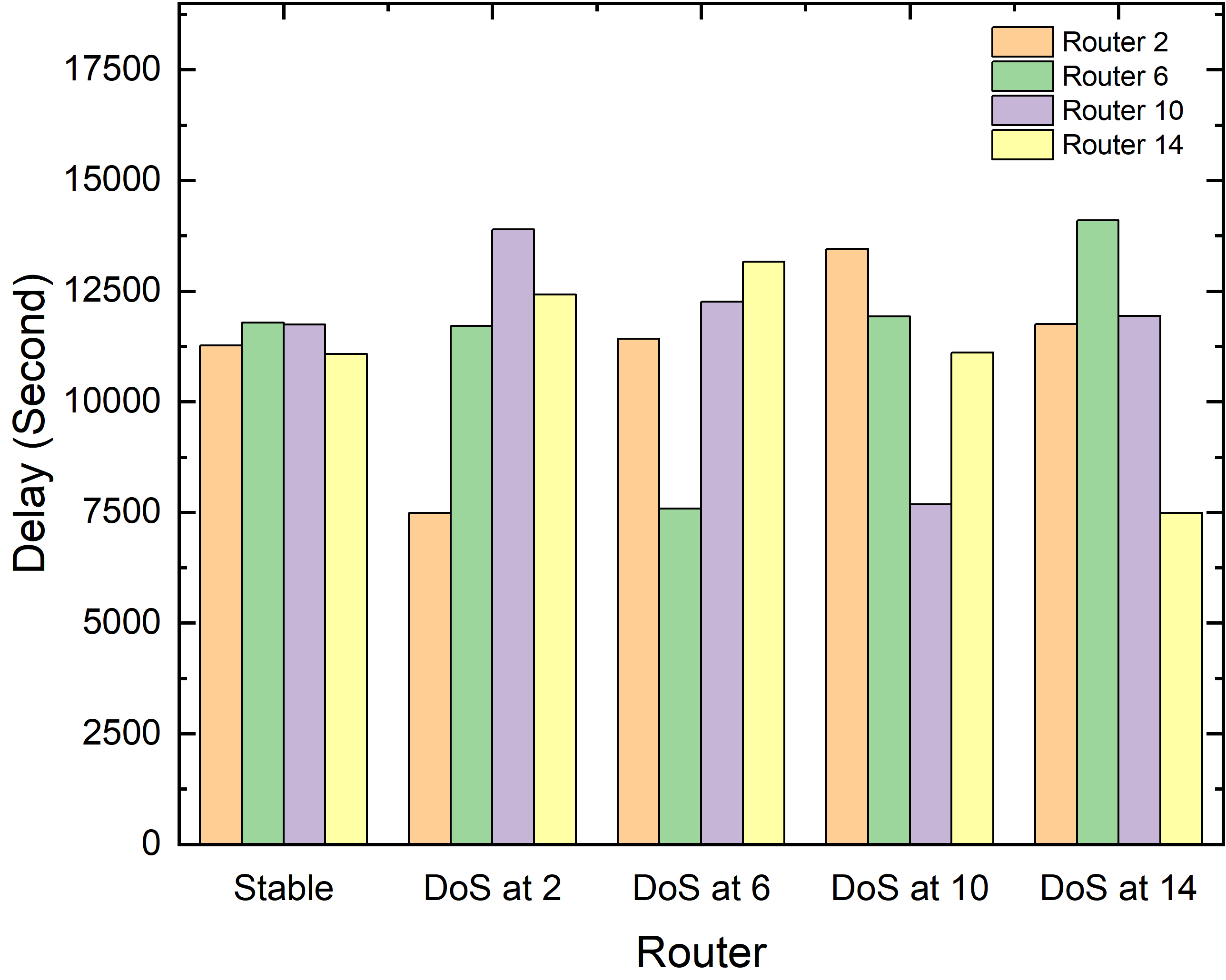}
    \caption{Case 3}
    \label{fig:fourth}
\end{subfigure}
        
\caption{Delays in the critical routers for each case with SimPy simulation. 
%the delay is set at an %arbitrary time unit, 
The `Stable' label denotes the base case (no attack) results for each case, and the horizontal axes showcase the different delays for each of the DoS/DDoS threats ran on each case (i.e., \textit{DoS at 5} in Figure \ref{fig:first} means a DoS threat on Router 5 in Case 1).}

\label{fig:delays}
\end{figure*}

The routers forward packets with a mean rate of 2.2 packets per second. A routing table is generated by creating branches to and from routers, with equal probabilities for each link. For example, a router connected to three links has a probability of 0.33 to forward the packet for each link. If a packet arrives in one of the links, that link does not forward it. When the packets arrive at the sink, they are destroyed. Each router is monitored during the simulation to measure how much time the router takes to forward a packet (i.e., the delay from the time a packet enters the router until the time the packet leaves the same router). This is called the router service time. The router records the service time of a packet every 0.5 seconds. The equation for delay has been given in Equation \ref{eq:delay}.
\begin{equation}
Delay = Queue~Waiting~Time / Number~of~Packets
\label{eq:delay}
\end{equation}
here $Queue~Waiting~Time$ represents the average waiting time for each packet. It is important to acknowledge a notable aspect observed in network simulations: the router directly connected to the sink typically exhibits a delay similar to that of the sink itself. This delay remains consistent and is unaffected by the network's topology, primarily dictated by the proximity between the router and the sink. In our simulation, routers calculate the delays of the packets they forward, rather than those they receive at the sink. Consequently, the router linked directly to the sink is not considered in the simulation.
To carry out the DoS attack, we rendered a router ineffective by drastically decreasing its probability of forwarding packets to an extremely minimal level. Normally, a router forwards packets with a certain probability, but by lowering this probability drastically, most packets are dropped instead of being forwarded. This mimics a DoS attack, causing the router to become a bottleneck and resulting in high packet loss and reduced network performance. In SimPy, after initializing the environment, the forwarding probability is set to 0.01. The environment is then run for a predefined duration to execute the attack.

\section{Results and Discussions}\label{Results}

In this section, we analyze and compare the different results of the three case studies. In the results, high network delay in a router means that a packet experiences high waiting times or service times as it travels through that router. Because a router forwards a packet in a first-in-first-out (FIFO) queuing process, high delay in a router means that many packets are arriving at that router and each packet has to wait their turn to be forwarded by the router. Based on the packet delays at the routers in the networks of cases 1, 2, and 3, the routers are ranked, and the top three routers are identified as critical elements since the ranking of critical elements in case 3 ends at three. Although this number may increase in the case of a large network, a more generalized analysis of the number of critical elements will be addressed in future work. %However, as 
Before the simulations, %as a first step, 
based on the topology of each network, we analyze the betweenness centrality, eccentricity centrality, eigenvector centrality and edge betweenness metrics for each router. 

\subsection{Critical Routers Based on Graph Theory Metrics}

Table \ref{tb2} shows the list of critical routers for each case. Based on the value of centrality, the critical routers of each case are noted – the higher the value the more critical the router except for eccentricity centrality where a low value indicates a critical router. 

\subsubsection{Case 1} In Case 1, using the cyber network of the IEEE 123-bus (Figure~\ref{fig: cyber network of IEEE 123}) with the centrality metric, router 7 is found to be the most critical one, followed by router 11 in second place. Similarly, with the help of edge betweenness centrality, the link between router (5, 7), router (7, 11) and router (10, 11) are found to be the most critical. From a network topology standpoint, these routers and links play a pivotal role in connecting various segments of the network. 

\subsubsection{Case 2} 
The central nodes, specifically routers 1 and 2 or 0, stand out as the most critical elements, as determined by the centrality metrics analysis of the radial topology (Figure \ref{fig:radialnetwork}). Following them are the next vital nodes, comprising routers 0 or 1, 2 or 3, 4, 5, and 6. Interestingly, nodes 3, 4, 5, and 6 lack redundancy. As a result, they demonstrate significance during the simulation phase.

\subsubsection{Case 3} 

In the ring topology illustrated in Figure \ref{fig:RingdNetwork}, routers 6, 10 are the critical routers in terms of betweenness centrality and eigenvector centrality. At the same time, routers 1, 2, 6, 10, and 14 (all routers) stand out as the most vital nodes in terms of eccentricity centrality.  Following closely are routers 2 and 14, comprising the next set of critical elements in terms of betweenness centrality and eigenvector centrality. Each of these routers, essential to the ring architecture, connects to three traffic generator nodes, except for router 1.

The next part of our experiment is to initiate a DoS attack on these critical routers. For the three use case network topologies, the delay packets experienced in each router are calculated after this event and plotted in Figure~\ref{fig:delays}. The results are described in the next section, and compared with the centrality results from Table~\ref{tb2}.

\subsection{Simulation of Use Case Results}

\subsubsection{Case 1} 

The results for use case 1 are shown in Figure~\ref{fig:first}. The routers that have larger delays are considered critical routers. Though the centrality metrics indicate the critical routers and links between the routers in a network, it fails to account for the position of those routers near the source and the sink of traffic. Assuming stable conditions, i.e., just normal traffic, the routers that had larger delays were Router 11 with the highest delay, followed by Router 5 and 9.  These routers experience the most delays because they are near the source, and the routers that are closer to the sink experience the least amount of delays. As a result, routers 2 and 6 which are not identified as critical routers by the betweenness centrality metrics also had significant delays and are included in our results. On the other hand, routers 4 and 8 showed the least amount of delays, and are removed from the plot. 

For our DoS disturbance simulation, routers 5, 9 and 11 are taken out independently since those are identified as the three most critical routers during the simulation in the stable case. In each instance, the router under DoS disturbance showed decreased delays (close to zero) since they are out of the system and no packets are going through them. For the DoS threat on Router 5, if we compare the router delays with the delays in stable conditions, it can be seen that the most affected router is router 2, which increased its delays by almost four times, with delays from less than 500 seconds in a stable network to more than 2000 seconds during the DoS disturbance on router 5. Note that router 2 is close to the source node G12. 
For DDoS, routers 5, 7 and 11 are taken out simultaneously from the system because the edge between 5-7 and 7-11 are identified as the critical edges from the edge betweenness centrality matrix. In this instance, routers 2, 9 and 10 are forwarding the bulk of the traffic since those routers are just one hop away from the generators G12, G16 and G17.

\subsubsection{Case 2} 
The results of the simulation from Figure~\ref{fig:third} (break marks have been added for better visualization of routers with smaller delays) show the highest delay values at routers 3, 4, 5, and 6. These routers are just two hops away from the generators. From the centrality matrix, these routers come in either second or third place in terms of criticality.  As a result, a DoS disturbance was initiated on router 3 and router 6 independently to predict the traffic patterns of the left-hand side and right-hand side of the topology respectively. In this instance, it is observed that for DoS disturbance the router which has been taken out from the system provides the least amount of delay and it is close to zero. When router 3 suffers a DoS disturbance the delay at router 6 is decreased indicating that more traffic is passing through the right-hand side of the network and the opposite result is observed when router 6 undergoes a DoS disturbance. Besides, the edges between routers 1, 3 and router 2, 6 come in second place on the edge betweenness centrality matrix. Thus, a DDoS disturbance was done at routers 1 and 3; followed by another DDoS at routers 2 and 6. In these instances, similar results are observed too when compared with the DoS disturbance. Since routers 1 and 3 are under disturbance and are situated on the left side of the network topology more packets are passed through the right side of the network. Again, when routers 2 and 6 are simultaneously taken out from the system, the left side of the network does the bulk transfer of packets.

\subsubsection{Case 3}

The results for case 3 are depicted in Figure~\ref{fig:fourth}. Under normal traffic conditions, all routers exhibit relatively similar delays. However, individual  DoS attacks were executed on routers 2, 6, 10, and 14 to evaluate their impact. When router 2 experienced a DoS attack, there was a decrease in delay in R2.
%, although a less significant decrease compared to cases 1 and 2.
This decrease in delay led to increased delays on other routers such as Routers 10 and 14. Similar results were observed with the disruption of routers 6, 10, and 14. Nonetheless, due to the small number of routers in this network and the fact that most of them are connected to sources of traffic, the decrease in delay in the disrupted routers of this use case are not as pronounced as in cases 1 and 2. DDoS disruption was not triggered in this scenario to avoid compromising the entire system, given the small number of routers this network has.

\section{Conclusions}\label{Conclusion}

In this paper, we employed SimPy to simulate traffic patterns for three use cases, with different communication topologies, and we
%used this simulation setup to 
%successfully 
%run 
ran DoS threats on different routers, 
%which were 
dictated by centrality metric analysis. 
%As such, our most essential contribution is that 
The
%is 
%work
analysis
%gave 
quantitatively compared and ranked
%results that prioritize 
devices and links in these use cases
%quantified the risk and analyzed 
based on the impact of a cybersecurity threat by assessing the criticality
%vulnerability 
of different components in multiple network topologies. %This improves our
These studies are
%is analysis and these types of comparisons are
%ese studies are 
%an 
important
%step 
for improving cyber-physical system
%the 
%understanding and 
modeling capabilities and
%, and they
%.
%of 
%the 
% behavior. %With that, we were able 
%These comparisons
%are an important step 
to
%toward
%Results show that these studies allow us to 
evaluate the fidelity of different
%the conducted 
%vulnerability 
risk assessment methods. Future work aims to scale up these assessments to generalize the number of critical elements that should be considered. Besides, time delay attacks that manipulate network traffic to introduce disruptions or deceive systems will be incorporated by simulating intentional delays with SimPy by assessing effects on communication latency, data synchronization, and overall system behavior.
%For future work, we aim to scale up our environment. 
It is crucial to acknowledge certain considerations in our study. 
%and assumptions made during the study. 
First, the cyber network is modeled based on the Network Layer (Layer 3) of the OSI model. However, it is important to note that the models are constructed based on existing cyber networks, which inherently come with their limitations. The size of the network considered in our study is moderate, but scalability is feasible due to the minimal overhead associated with the simulation process. Besides, centrality metrics used do not differentiate well between node or edge types, such as source, sink, or router. Therefore, it is essential to recognize that our analyses consider nuance in this aspect. 

In summary, the combined use of graph theory and simulation experiments yields valuable insights into the security of cyber network topologies. However, to enhance the evaluation of node significance in communication networks, it is recommended to consider node types. Thus, employing multiple graph metrics can offer a more comprehensive understanding of network dynamics and vulnerabilities.

\bibliographystyle{IEEEtran}
\bibliography{ref}

\vfill

\end{document}